\documentstyle{article}

 \def\vt{t\kern-0.22em\raise.18ex\hbox{\char'47}\lower.18ex\hbox{}\kern-0.08em}
 
\newtheorem{th}{Theorem}[section]
\newtheorem{co}{Corollary}[section]
\newtheorem{lm}{Lemma}[section]

\newcommand{\old}[1]{{}} 

\newcounter{obr}
\newcounter{tabul}
\begin{document}
\title{ A new  upper bound  for  the   clique cover number  with applications
\thanks{This paper appeared in Congressus Numerantium, 205 (2010), 105-111 .} 
}
\author{Farhad Shahrokhi\\
Department of Computer Science and Engineering,  UNT\\
farhad@cs.unt.edu
}

\date{}
\maketitle
\thispagestyle{empty}
\date{} \maketitle


\begin{abstract}
Let $\alpha(G)$ and  $\beta(G)$, denote  the size  of a largest 
independent set  and the clique cover number of an undirected graph  $G$. 
Let   $H$ be an interval graph with $V(G)=V(H)$ and $E(G)\subseteq E(H)$, 
and let $\phi(G,H)$   
denote the maximum of ${\beta(G[W])\over \alpha(G[W])}$
overall induced subgraphs $G[W]$  of $G$  that are  
cliques in $H$. 
The main result of this paper is to prove that for any graph $G$  
$${\beta(G)}\le 2 \alpha(H)\phi(G,H)(\log \alpha(H)+1),$$
where, $\alpha(H)$ is the size of a largest independent
set in  $H$.
We further provide a generalization 
that significantly  unifies or  improves 
 some  past  algorithmic and structural  results concerning 
the clique cover number for some  well known  intersection  graphs.

\end{abstract}

\section{Introduction}
Throughout this paper  $G=(V(G),E(G))$ is a simple undirected graph  
with $|V(G)|=n$. 
Let  $H$ be a subgraph  of $G$. We denote  by $\beta(H)$   the 
clique cover number of $G$. 
We further denote by  $\alpha(H)$ and $\omega(H)$  the cardinality of 
a largest independent set and a largest clique in  $H$, respectively. 
Let $G[W]$ denote the induced graph in $G$  on the vertex set  $W\subseteq V$.

Clearly $\beta(G)\ge \alpha(G)$, for any $G$. In addition,  
there exit graphs $G$  with arbitrary large  $\beta(G)$ so that  
$\alpha(G)\le 2$ \cite{My}. 
Consequently, given a graph $G$,  one can not expect 
to  have an upper bound for $\beta(G)$  that involves
$\alpha(G)$ only.  Furthermore,
the problems of  computing    $\beta(G)$ exactly for a graph $G$, 
or even approximating it,    
are known to be computationally hard \cite{GJ}, \cite{BO}.
Hence, deriving   any general upper  bound for the clique cover
number  that may offer algorithmic
consequences is a  significant move forward.

Throughout this paper we say a graph $H$ is a  {\it supergraph} of $G$, if
$V(G)=V(H)$ and  $E(G)\subseteq E(H)$. 
 Let $H$ be a supergraph of $G$. Define $\phi(G,H)$ as 

\[\phi(G,H)=\max_{{W\subseteq V\atop H[W]~is~a~clique}}\Bigg\{{\beta(G[W])\over\alpha(G[W])}\Bigg\}.\]

An interval graph is the intersection graph of a  set of intervals on the 
real line \cite{To},\cite{co}. 
The main contribution this paper is to prove the  following. 

\begin{th}\label{t1}
{\sl  
Let  $H$  be  an interval  
supergraph of  $G$, then 
\[\beta(G)\le 2\alpha(G)\phi(G,H)(\log(\alpha(H))+1).\]

}
\end{th}
Note that $\alpha(H)\le\alpha(G)$ and hence Theorem 1.1 implies that
$\beta(G)\le 2\alpha(G)\phi(G,H)(\log(\alpha(G))+1)$, for any graph
$G$.
\vskip .3cm

We further   generalize  Theorem \ref{t1}  to the case 
$E(G)=\cap_{i=1}^t E(H_i)$, where  the supergraph $H_i,i=1,2,...,t-1,$   
is an interval graph and the supergraph $H_t$ is a perfect graph. 
 This generalization 
significantly  unifies or  improves  
some   past  algorithmic and structural  results concerning 
the clique cover and the maximum independent set problems  
for some  well known  intersection  graphs. 
\vskip .3cm

Specific applications include   generalizing the results  
on the traversal number of rectangles \cite{Ka}, \cite{Chan}, 
\cite{N}, drastically 
improving an upper bound  for the  clique cover number  of 
interval filament graphs\cite{CC} and improving the approximation 
factor  for the polynomial time  approximation of the clique cover number of 
interval filament graphs \cite{KS}.
\vskip .3cm

The strength of the work presented here is the generality that leads
to the   significant unification of several past results.
In addition, the  underlying method   utilizes     
  simple  properties associated with the  linear orderings 
 and the clique separation  of  interval graphs,  and thus  may be 
of an  independent interest.

\section{Main Result}
Two disjoint subsets $A,B\subseteq V(G)$
are {\it separated} in $G$, if there are no edges $ab\in E(G)$ with $a\in A$ 
and $b\in B$.  
Let $\pi:v_1,v_2,...,v_n$ be  linear ordering  of $V(G)$.
 For $i=1,2,...,n$ the sets  
$\{v_1,v_2,...,v_i\}$ and $\{v_i,v_{i+1},...,v_n\}$ are  denoted by 
$V_i$ and $V^i$, respectively. 
We define $V_0=\emptyset$. 
For $i=1,2,...,n$, let  
$N(v_i)$ denote the set of all vertices adjacent to $v_i$ in $V_{i-1}$.

Let $G$ be an interval graph whose interval representation is   
$J=\{[a_i,b_i]|i=1,2...,n\}$. 
Thus $G$ is the intersection graph of intervals in $J$. 
Let $\pi$ be a linear  ordering of elements in $J$ in the increasing
order of $a_i$, for $i=1,2,...,n$. Then,  $\pi$ is called
a  {\it canonical linear ordering}.
The following Lemma  states some  known  properties  concerning the canonical
linear orderings of interval graphs \cite{To}, \cite{co}.

\begin{lm}\label{l1}
{\sl Let  $\pi:v_1,v_2,...,v_n$ be a canonical linear  ordering  of vertices 
of an interval graph $G$. Then, 
for any $i=1,2,...,n$, 
$N(v_i)$  is a clique in $G$. 
In addition,  
for any $i=1,2,...,n$, 
$V_{i-1}-N(v_i)$
is separated from $V^i$ in $G$.  
}
\end{lm}

We  now prove  the main result.

{\bf Proof of Theorem 1.1} 
We will show that there is a clique cover $C(G)$, and an independent
set $I(G)$ so that 
\[|C(G)|\le 2|I(G))|\phi(G,H)(\log(\alpha(H))+1). \]

To prove the claim,  we induct on $\alpha(H)$. 
Our proof gives rise to a recursive divide and conquer algorithm
for constructing $C(G)$ and $I(G)$. 
For $\alpha(H)=1$, the claim is a direct consequence of the definitions, since in this case
$\beta(H)=\alpha(H)=1$.  Thus  $H$ is a complete graph, and consequently
$\beta(G)\le \alpha(G)\phi(G,H)$.
Now let $k\ge 2$, assume that the claim is valid for all graphs 
$G$ and their interval supergraphs $H$ with $\alpha(H)<k$, and let   
$H$ be an interval supergraph of $G$ with $\alpha(H)=k$. 
Let $\pi$ be a canonical linear ordering, 
let $I, |I|=\alpha(H)$ be a maximum independent 
set in $H$, let $j$ be the index of  $\lfloor{\alpha(H)\over 2}\rfloor$-th  
vertex in   $I$ in the ordering  induced by $\pi$. Finally, let
$j^*$  be the index of the vertex in $I$ that appears immediately
after   $v_j$,   in the linear ordering induced by $\pi$.     
Let  $V'_j=V_j-N(v_{j^*})$ and note that $H[V'_j]$ and $H[V^{j^*}]$ are interval
graphs that contain  independent  sets of sizes 
$\lfloor{\alpha(H)\over 2}\rfloor$ and 
$\alpha(H)-\lfloor{\alpha(H)\over 2}\rfloor=\lceil{\alpha(H)\over 2}\rceil$, respectively.
It follows that $\alpha(H[V'_j])=\lfloor{|\alpha(H)|\over 2} \rfloor$, and   
$\alpha(H[V^{j^*}])=\lceil{\alpha(H)\over 2}\rceil$, since by 
Lemma \ref{l1} $V'_j$ is separated  from $V^{j^*}$ in $H$. 
Note further that   $H[V'_j]$ and  $H[V^{j^*}]$ are   supergraphs of  
$G[V'_j]$ and $G[V^{j^*}]$, respectively, and hence
$V'_j$ and $V^{j^*}$ are also  separated in $G$. In addition,   
$\phi(G[V'_j],H[V'_{j}])\le\phi(G,H)$ 
and $\phi(G[V^{j^*}],H[V^{j^*}])\le\phi(G,H)$. 
Let $C(G[V'_j])$ and $I(G[V'_j])$ be clique cover and independent sets that are
obtained by the application of the   induction hypothesis to $G[V'_j]$ and $H[V'_j]$, then
$$|C(G[V'_j])|\le 2|I(G[V'_j])|\phi(G,H)(\log\lfloor{\alpha(H)\over 2}\rfloor+1).$$
Similarly, let $C(G[V^{j^*}])$ and $I(G[V^{j^*}])$  be the clique  cover and 
the independent set obtained   by applying the induction hypothesis to $G[V^{j^*}]$ 
and $H[V^{j^*}]$. Then,
$$|C(G[V^{j^*}])|\le 2|I(G[V^{j^*}])|\phi(G,H)(\log\lceil{\alpha(H)\over 2}\rceil+1).$$
Next, note that   $H[N(v_{j^*})]$ is a clique in $H$, and thus, there is a 
clique cover $C(G[N(v_{j^*})])$ and an independent set $I(G[N(V_{j^*})])$ so
that   $|C(G[N(v_{j^*})])|\le\phi(G,H)|I(G[N(v_{j^*})])|$. Now, let  
$C(G)=C(G[V^{j^*}])\cup C(G[N(v_{j^*})])\cup C(G[V'_j])$ and   observe that  
$|C(G)|=|C(G[V'_j])|+|C(G[V^{j^*}])|+|C(G[N(v_{j^*})])|$,
and hence by combining the last three inequalities, we obtain 
{\small \[|C(G)|\le 2\phi(G,H)(|I(G[V_j])|+|I(G[V^{j^*}])|)(\log\lceil{\alpha(H)\over 2}\rceil+1)+\phi(G,H)I(G[N(v_{j^*})]).\]}
Next  define  $I(G)$ to be the larger of the two    independent sets   
$I(G[V^{j^*}])\cup I(G[V'_j])$,  $I(G[N(v_{j^*})])$. Then, the last inequality gives
$$|C(G)|\le 2\phi(G,H) |I(G)|(\log\lceil{\alpha(H)\over 2}\rceil)+{3\over 2}.$$
To finish the proof, observe that
$\log (\lceil{i\over 2}\rceil+{3\over 2})\le 
\log(i)+{1}$ for 
any integer  $i,i\ge 2$. 

$\Box$

An immediate consequence of Theorem    \ref{t1}  is the following.
\begin{co}\label{co1}
{\sl 
Let $H_1$  and $H_2$ be 
supergraphs  of $G$, that are an interval graph   and a  perfect
graph, respectively.
If $E(G)=E(H_1)\cap E(H_2)$, then  there is a  clique cover
$C(G)$ and an independent set $I(G)$, in $G$,  so that  
\[|C(G)|\le 2|I(G)|(\log(\alpha(H_1))+1).\] 

Moreover, $I(G)$ and $C(G)$ can be computed in polynomial time.
}
\end{co}
{\bf Proof.} Let $W\subseteq V(G)$
so that  $H_1[W]$ is a clique.  Then $G[W]=H_2[W]$, since $E(G)=E(H_1)\cap E(H_2)$.   Thus 
 $\beta(G[W])=\alpha(G[W])$, since  $H_2$ 
is a perfect graph.  
This implies that  $\phi(G,H_1)=1$, and the claims for 
the upper bound on $C(G)$ follow.  To finish the proof,  observe that 
for any subgraph of a perfect graph the clique cover and the 
independent set can be computed in polynomial time \cite{GLS}. 
$\Box$
\vskip .3cm

{\bf Remark.} The time complexity of constructing $I(G)$ and $C(G)$ is
$O({\hat T}(n)+n^2)\log(n)),$ where ${\hat T}(n)$ is the time
that it takes to  compute   the clique cover number and the maximum 
independent  set in any subgraph of the perfect graph $H_2$. 
Note that  $\hat T(n)$ is a polynomial of $n$, for the  
perfect graph $H_2$. 
However, the  degree of the
polynomial depends on the structure of the $H_2$. For instance,
$\hat T(n)=O(n^2)$, when $H_2$ is an incomparability graph.

Using induction one can generalize Corollary  \ref{co1}.
\begin{co}\label{co2}
{\sl 
Let $t\ge 2$, let $H_1,H_2,...,H_{t-1}$ be  interval
super graphs  of   $G$, and let $H_t$ be a perfect supergraph of $G$.
If $E(G)=\cap_{i=1}^t E(H_i)$,  then there is a vertex  cover $C(G)$ 
and independent set $I(G)$, in $G$, so that

\[|C(G)|\le 2^{t-1}|I(G)|{\Pi_{i=1}^{t-1}(\log(\alpha(H_i))+1)}.\]
Moreover, $I(G)$ and $C(G)$ can be computed in polynomial time. 
}
\end{co}

\section{Applications}
In this section we point out  a few applications.
Gavril \cite{gav}  introduced   the class of 
{\it interval filament} graphs
which is an important class of intersection graphs and has shown
that for any graph $G$ in this class $E(G)=E(H_1)\cap E(H_2)$, where
$H_1$ and $H_2$ are supergraphs of $G$ which are an interval graph and
an incomparability graph, respectively. 
This important  class of graphs contains
the classes  of incomparability, polygon circle, chordal, circle, circular arc,
and outer planar graphs. Gavril \cite{gav} has shown that $\alpha(G)$ and $\omega(G)$ can be computed in low-order polynomial time for   
 any  interval filament graph $G$.   
However, since  the problem of computing 
$\beta(G)$ for a circle graph $G$ is 
$NP-hard$ \cite{KS}, computing the clique cover number  is also $NP-hard$ 
for interval filament graphs. 

Cameron and Hoang \cite{CC} have shown that 
for any interval filament graph $G$,  $\beta(G)=O(\alpha(G)^2)$. 
One can also obtain this  result by the application of the general  
method of    Pach and    T\"or\"ocsik \cite{pt}.
The application of  Corollary \ref{co1} to an interval filament graph $G$ 
gives $\beta(G)=O(\alpha(G)\log\alpha(G))$, 
a  significant  improvement over  the upper bound
in \cite{CC}. Since  our method  in Corollary \ref{co1}  constructs  
a  suitable clique  cover $C(G)$  in polynomial time,
it also improves the best known  approximation factor 
for the  polynomial time  approximation of $\beta(G)$ in  \cite{KS} 
from $O(\log(n))$ to  $O(\log(\alpha(G)))$.  
We remark that the $O(\log(n))$ approximation factor of  Keil and Stewart 
in \cite{KS} has been
derived for a larger class of graphs, namely  the class of  subtree filament graphs,  that contains the class of  interval filament graphs.
\vskip .3cm

Kostochka and Kratochvil \cite{KK} have shown 
$\beta(G)=O(\alpha(G)\log(\alpha(G)))$ for any polygon circle graph $G$. 
Since any polygon circle graph is also
an interval filament, Corollary \ref{co1} also implies the result in \cite{KK}.   
\vskip .3cm

Let $G$ be the intersection graph of a set $S$  of axis parallel rectangles 
in the  plane.
Karolyi \cite{Ka} was the first to 
prove that
$\beta(G)=O(\log(\alpha(G))\alpha(G)).$ (For a simpler proof see \cite{FK}.)
Later, Agarwal et al\cite{AG} studied the related  problem of computing 
$\alpha(G)$ which is known to be $NP-hard$,  and  designed an $O(n\log(n))$ 
algorithm to construct  a  an upper  bound of $O(\alpha(G)\log(n))$. Different  algorithmic variations 
of the method of   Agarwal et al have been  studied in   
\cite{Chan},  particularly in view  
of the applications of computing $\alpha(G)$ in  map labeling and cartography. 
It  is readily seen that $E(G)=E(H_1)\cap E(H_2)$, where  
$H_1$ and $H_2$ are interval supergraphs of $G$. Specifically,
the interval  orders  $\prec^1$ and $\prec^2$  associated with $H_1$    
and $H_2$ define  the separation properties of rectangles
in $S$, in the horizontal and vertical directions. Since any interval graph
is perfect,  Corollary \ref{co1} applies and gives 
$\beta(G)=O(\alpha(G)\log(\alpha(G)))$. 
In addition, Corollary  \ref{co2} gives
$\beta(G)=O(\alpha(G)\log^t(\alpha(G)))$ and  constructs a  suitable 
clique cover and an independent set in polynomial time,  
which is the same as the best known result in \cite{N},  
for the generalization of the problem to $R^t$, $t\ge 2$. 
A similar result follows from Corollary \ref{co2},
when  $G$ is  the intersection graph of a set $S$ of   convex polygons in the 
the plane,  where each  $P\in S$ is obtained  by the translation
and  magnification of a particular convex polygon with $t$ corners, 
for a constant $t\ge 3$.


\begin{thebibliography}{99}


\bibitem{AG}
Agarwal P.K., Kreveld M., Suri S., Label placement by maximum independent sets in rectangles, 
{\it Comput. Geometry:Theory and Appl.}
 11(1998) 209-218.

\bibitem {BO}
Bellare M.,
 Goldreich O., Sudan M., Free bits, PCPs and non-approximability -
 towards tight results, SIAM Journal on Computing, 27(1998), 804-915.

\bibitem{Chan} 
Chan T., A note on maximum independent sets in rectangle intersection graphs, 
{\it Inform. Process. Lett.} 89 (2004) 19-23. 

\bibitem{co} Golumbic  M. C, Algorithmic Graph Theory and Perfect Graphs (Annals of Discrete Mathematics, Vol 57, North-Holland Publishing Co., Amsterdam, The Netherlands, 2004. 

\bibitem{CC}  Cameron K., Hoang C., On the  structure of certain intersection
graphs, IPL,  99 (2006) 59-63.


\bibitem{FK}
Fon-Der-Flaass D.G, Kostochka, A. V., Covering  boxes by points,
{\it Discr. Math.} 120(1993) 269-275.

\bibitem{gav} Gavril F, , Maximum weight independent sets and cliques in intersection graphs of filaments, IPL  73 (2000)  181-188.

 \bibitem{GLS} Gr\"otsche M. , Lovas  L. , Schrijver A., The ellipsoid method and its consequences in combinatorial optimization, Combinatorica 1 (1981) 169-197. 


\bibitem{KS} Keil M.J, Stewart L., Approximating the minimum clique 
cover and other hard problems in subtree filament graphs,  {\it Discrete Applied Mathematics}, 154 (2006) 1983-1995.

\bibitem{KK}
 Kratochvil J.,  Kostochka
 A., Covering and coloring polygon circle graphs,
{\it Disc. Math.} 163 (1997) 299-305.

\bibitem{Ka}
Karolyi G., On point covers of parallel rectangles, {\it Periodica Mathematica 
Hungarica} 23 (1991) 105-107.
\bibitem{GJ}
Garey M. R., Johnson D. S., {\it Computers and intractability: a guide to NP-Completeness}, W. H. Freeman and Co., 1978. 
\bibitem{N}
Nielsen F., Fast stabbing of boxes in high dimensions, TCS, 246 (2000)
53-72. 

\bibitem{My}
Mycielski, J.  Sur le coloriage des graphes, Colloq. Math. 3 (1995) 161–162.

\bibitem{pt}
 Pach J. and  T\"or\"ocsik J., 
Some geometric applications of Dilworth's theorem, 
{\it Disc. Comput. Geometry}  21(1994), 1--7.

\bibitem{To}  Trotter W.T.,   New perspectives on interval orders and interval
graphs,  in Surveys in Combinatorics, Cambridge Univ. Press (1977) 237-286.

\bibitem{To2}  Trotter W.T., Combinatorics and partially ordered sets: Dimension theory, Johns Hopkins series in the mathematical sciences, The Johns Hopkins University Press, 1992.



































\end{thebibliography}
\end{document}